\documentclass[a4,12pt,amsmath]{article}
\usepackage[headings]{fullpage}
\usepackage{amssymb,epic,eepic,epsfig,amsbsy,amsmath}
\numberwithin{equation}{section}
\usepackage{amscd}
\usepackage{times}
\usepackage{amssymb}
\def\SO{\mathop{\mathrm {SO}}\nolimits}
\def\SU{\mathop{\mathrm {SU}}\nolimits}
\def\HP{\mathop{\mathrm {HP}}\nolimits}
\def\K{\mathop {\mathrm {K}}\nolimits}
\def\E{\mathop {\mathrm {E}}\nolimits}

\def\Cyl{\mathop{\mathrm {Cyl}}\nolimits}
\def\Cone{\mathop{\mathrm {Cone}}\nolimits}
\def\ev{\mathop{\mathrm {ev}}\nolimits}
\def\S{\mathop{\mathbf {S}}\nolimits}
\def\b{\mathop{\mathbf {b}}\nolimits}
\def\t{\mathop{\mathbf {t}}\nolimits}
\def\N{\mathop{\mathbf {N}}\nolimits}

\def\Tot{\mathop{\mathrm {Tot}}\nolimits}
\def\1{\mathop{\mathbf {1}}\nolimits}

\def\H{\mathop{\mathrm {H}}\nolimits}
\def\d{\mathop{\mathbf {d}}\nolimits}
\newtheorem{thm}{Theorem}[section]
\newtheorem{defn}[thm]{Definition}

\newtheorem{cor}[thm]{Corollary}
\newtheorem{lem}[thm]{Lemma}
\newtheorem{rem}[thm]{Remark}
\em
\begin{document}
\title{Category of Noncommutative CW Complexes. II\footnote{\bf Version  from \today .
}\footnote{\bf The work was supported in part by Vietnam National
Project for Research in Fundamental Sciences and was completed
during the stay in June - July, 2007 at the Abdus Slam ICTP,
Trieste, Italy }}
\author{Do Ngoc Diep}
\maketitle
\begin{abstract}
We introduce in this paper the notion of noncommutative  Serre
fibration (shortly, NCSF) and show that up to homotopy, every NCCW
complex morphism is some noncommutative Serre fibration. We then
deduce a six-term exact sequence for the periodic cyclic homology
and for K-theory of an arbitrary noncommutative Serre fibration.
We also show how to use this technique to compute K-theory and
cyclic theory of some noncommutative quotients.
\end{abstract}
\section*{Introduction}
It is well-known that for short exact sequences of closed
two-sided ideals and C*-algebras of type
\begin{equation} \CD 0
@>>> A @>>> B @>>> C @>>> 0\endCD\end{equation}
 where $A$ is some
(closed) two-sided ideal in $B$ and $C \cong B/A$ the quotient
C*-algebra, the periodic cyclic homology admits six-term exact
sequences of type
\begin{equation}\CD \HP_*(A) @>>> \HP_*(B) @>>> \HP_*(C)\\
@A\partial_{*+1} AA @. @VV\partial_* V\\
\HP_{*+1}(C) @<<< \HP_{*+1}(B) @<<< \HP_{*+1}(A)
\endCD\end{equation}
The same is true for K-theory
\begin{equation}\CD \K_*(A) @>>> \K_*(B) @>>> \K_*(C)\\
@A\partial_{*+1} AA @. @VV\partial_* V\\
\K_{*+1}(C) @<<< \K_{*+1}(B) @<<< \K_{*+1}(A)
\endCD\end{equation}
It is natural to ask whether the condition of $A$ being an ideal
is necessary for this exact sequence. We observe that this
condition is indeed not necessarily to be satisfied. It is enough
to consider the relative generalized homology in place of the
generalized homology of the quotient algebra, what we don't have
in general for a pair of algebra and subalgebra. We use the
corresponding notion of noncommutative mapping cylinder,
noncommutative mapping cone and suspension etc. introduced in
\cite{diep1} in order to manage the situation. The ideas are
therefore originated from algebraic topology.

For more special pairs, some thing like ordinary Serre fibrations
with the {\it homotopy lifting property} (HLP), we have more
properties. For this we arrive to a more general condition of a
{\it noncommutative Serre fibration} (NCSF) as some homomorphism
of algebras with the so called HLP (Homotopy Lifting Property).
But we do restrict to consider only the so called {\it
noncommutative CW- complexes}, \cite{diep1}. We prove that up to
homotopy we can change any homomorphism between noncommutative
CW-complexes in order to have a noncommutative Serre fibration.
This main theorem let us then to apply it to study noncommutative
decomposition series for operator algebras, started in
\cite{diep2}

Let us describe the contents of the paper. In Section 1 we define
some noncommutative object like noncommutative (NC) mapping
cylinder, NC mapping cone, NC suspension, etc....  and the HLP
(Homotopy Lifting Property) and NCCW complex. As the main result,
we prove that in the category of NCCW complexes, every morphism is
homotopic to a Serre fibration. We deduce in Section 2 the hexagon
diagram for the periodic cyclic homology and K-theories of NC
Serre fibration. In the ordinary algebraic topology one uses
namely some computation of cohomology of spheres to compute
homology of orthogonal and unitary groups. In our noncommutative
theory, the cyclic homology of quantum orthogonal and quantum
unitary groups are known from the representation theory we then
use those to compute homology of quantum spheres as noncommutative
quotients. In a next paper we prove also that for NC Serre
fibrations there is also spectral sequences converging to HP and K
theories with $E_2$ term like $\E^2_{p,q} =\HP_p(B,A; \HP_q(A))
\Longrightarrow \HP_{p+q}(B)$, and $\E^2_{p,q} = \HP_p(B,A;K_q(A))
\Longrightarrow \K_{p+q}(B)$. In the works
\cite{diep2}\cite{noncomm} we reduced the problem of studying the
structure of an arbitrary GCR C*-algebra (of type I) to some
filtration by a towers of ideals with good enough quotients, what
permit to use spectral sequence ideas into studying. The paper is
the revised version of the preliminary version in
arXiv:math.QA/0211048.

\section{Noncommutative Serre Fibrations}
We start with the notion of NCCW, introduced by S. Eilers, T. A.
Loring and G. K. Pedersen \cite{elp} and G. Pedersen
\cite{pedersen}, see also \cite{diep1}.

\begin{defn}\label{defn1.2} A dimension 0 NCCW complex {\rm is defined, following
\cite{pedersen} as a finite sum of C* algebras of finite linear
dimension, i.e. a sum of finite dimensional matrix algebras,}
\begin{equation}A_0 = \bigoplus_{k} \mathbf M_{n(k)}.\end{equation}
In dimension n, an NCCW complex  {\rm is defined
 as a sequence $\{A_0,A_1,\dots,A_n\}$ of C*-algebras $A_k$
 obtained each from the previous one by the pullback construction
 \begin{equation}\CD 0 @>>> I^k_0F_k @>>> A_k @>\pi>> A_{k-1} @>>> 0\\
 @. @| @VV\rho_k V  @VV\sigma_k V @. \\
 0 @>>> I_0^kF_k @>>> I^kF_k @>\partial>> \mathbf S^{k-1}F_k @>>>
 0,\endCD \end{equation}
 where $F_k$ is some C*-algebra of finite
 linear dimension, $\partial$ the restriction morphism, $\sigma_k$
 the connecting morphism, $\rho_k$ the projection on the first
 coordinaes and $\pi$ the projection on the second coordinates in
 the presentation
 \begin{equation}
 A_k = \mathbf I^kF_k \bigoplus_{\mathbf S^{k-1}F_k}A_{k-1}\end{equation}
 }\end{defn}

\begin{defn}{\rm
We say that the morphism $f: A \longrightarrow B$ admits the so
called {\it Homotopy Lifting Property (HLP)} if for every algebra
$C$ and every morphism $\varphi : A \longrightarrow C$ such that
there is some morphism $\tilde{\varphi} : B \longrightarrow C$
satisfying $\varphi = \tilde{\varphi}\circ f$, and for every
homotopy $\varphi_t : A \longrightarrow C$, $\varphi_0 = \varphi$,
there exists a homotopy $\tilde{\varphi}_t : B \longrightarrow C$,
$\tilde\varphi_0 = \tilde\varphi$,  such that for every $t$,
$\varphi_t = \tilde{\varphi}_t \circ f$, i.e. the following
diagram is commutative }
\begin{center}\begin{equation}
\begin{picture}(60,60)
\put(0,0){A} \put(50,0){C} \put(0,50){B}
\put(10,5){\vector(1,0){35}} \put(5,10){\vector(0,1){35}}
\put(15,45){\vector(1,-1){30}} \put(10,25){f}
\put(37,25){$\tilde{\varphi}_t,\tilde{\varphi}_0=\tilde{\varphi}$}
\put(25,10){$\varphi_t$} \put(5,-10){$\varphi_0 =
\tilde{\varphi}_0 \circ f$}
\end{picture}
\end{equation}\end{center}
\end{defn}
\begin{defn}{\rm
A morphism of C*-algebras $f : A \longrightarrow B$ with HLP axiom
is called a {\it noncommutative Serre fibration} (NCSF).}
\end{defn}

\begin{thm}\label{thm1}
In the category of NCCW complexes, for every morphism $f : A
\longrightarrow B$, there is some homotopies $A \sim A'$, $B \sim
B'$ and a morphism $f' : A' \longrightarrow B'$ which is a NC
Serre fibration.
\end{thm}

Before prove this theorem we do recall \cite{diep1} the following
notions of noncommutative cylinder and noncommutative mapping
cone.

\begin{defn}[NC cone]{\rm For C*-algebras the} NC cone {\rm of $A$
is defined as the tensor product with } $\mathbf C_0((0,1])$, i.e.
\begin{equation}\Cone(A) := \mathbf C_0((0,1]) \otimes
A.\end{equation}
\end{defn}

\begin{defn}[NC suspension]{\rm For C*-algebras the} NC suspension
{\rm of $A$ is defined as the tensor product with } $\mathbf
C_0((0,1))$, i.e. \begin{equation}\mathbf S(A) := \mathbf
C_0((0,1))\otimes A.
\end{equation}
\end{defn}

\begin{rem}
If $A$ admits a NCCW complex structure, the same have the cone
$\Cone(A)$ of $A$ and the suspension $\mathbf S(A)$ of $A$.
\end{rem}

\begin{defn}[NC mapping cylinder]{\rm
Consider a map $f : A \to B$ between C*-algebras. The NC {\it
mapping cylinder} $\Cyl(f:A\to B)$ is defined by the pullback
diagram \cite{diep1}
\begin{equation}\CD \Cyl(f) @>pr_1 >> \mathbf C[0,1] \otimes A \\
@V{pr_2}VV  @VV{f\circ \ev(1)}V\\
B @>id>> B\endCD,\end{equation} where $\ev(1)$ is the map of
evaluation at the point $1\in [0,1]$. It can be also defined
directly as follows.
 In the algebra $\mathbf C(\mathbf I) \otimes A
\oplus B$ consider the closed two-sided ideal $\langle \{1\}
\otimes a - f(a), \forall a\in A\rangle $, generated by elements
of type $\{1\} \otimes a - f(a), \forall a\in A$. The quotient
algebra
\begin{equation}\Cyl(f)=\Cyl(f:A\to B) := \left(\mathbf C(\mathbf
I) \otimes A \oplus B\right)/ \langle \{1\} \otimes a - f(a),
\forall a\in A\rangle \end{equation} is called the} NC mapping
cylinder {\rm and denote it by } $\Cyl(f:A\to B)$.
\end{defn}

\begin{rem}
It is easy to show that $A$ is included in $\Cyl(f: A \to B)$ as
$\mathbf C\{0\} \otimes A \subset \Cyl(f:A\to B)$ and $B$ is
included in also $B \subset \Cyl(f:A\to B)$.
\end{rem}
\begin{defn}[NC mapping cone]{\rm The NC mapping cone
$\Cone(\varphi)$ is defined from the pullback diagram
\begin{equation}  \CD \Cone(\varphi) @>pr_1 >> \mathbf C_0(0,1] \otimes A \\
@V{pr_2}VV  @VV{f\circ \ev(1)}V\\
B @>id>> B\endCD,  \end{equation} where $\ev(1)$ is the map of
evaluation at the point $1\in [0,1]$. It can be also directly
defined as follows. In the algebra $\mathbf C((0,1]) \otimes A
\oplus B$ consider the closed two-sided ideal $\langle \{1\}
\otimes a - f(a), \forall a\in A\rangle $, generated by elements
of type $\{1\} \otimes a - f(a), \forall a\in A$. We define the}
mapping cone {\rm as the quotient algebra }
\begin{equation}\Cone(f)=\Cone(f:A\to B) := \left(\mathbf
C_0((0,1]) \otimes A \oplus B\right)/ \langle \{1\} \otimes a -
f(a), \forall a\in A\rangle .\end{equation}
\end{defn}
\begin{rem}
It is easy to show that $B$ is included in $\Cone(f: A \to B)$.
\end{rem}

\noindent {\sc Proof of Theorem \ref{thm1}.}

{\sc Step 1}. From the definition $\Cyl(f:A\to B)$, we see that
\begin{lem}
$B' = \Cyl(f:A\to B)$ is homotopic to B and $A \cong A \otimes
\mathbb C(\{0\}) \hookrightarrow B'.$ \end{lem} We can therefore
from now on suppose that $A \hookrightarrow B$.

{\sc Step 2}. Let us denote $\Pi(A,B)$ the space of all piece-wise
linear curve starting from $A$ and ending in $B$. Every element in
this space is la linear combination of elements from $B$, which is
a NCCW complex and therefore
\begin{lem} The algebra  $\Pi(A,B)$
 is also a NCCW complex in sense of Definition \ref{defn1.2}.
\end{lem}

{\sc Step 3}. The following lemma is clear.
\begin{lem}
$A':= \Pi(A,B)$ is homotopic to $A$.
\end{lem}
Indeed, by change of parametrization $\gamma(t) \mapsto
\gamma(st), t\in \mathbf I = [0,1]$, for $s\in \mathbf I = [0,1]$
we can produce a homotopy between curves and their starting
points.

 {\sc Step 4}. Let us consider the morphism $A' =
\Pi(A,B)\longrightarrow B$ which corresponds to each piece-wise
linear path its end-point.
\begin{lem} This map $A' \longrightarrow B$ satisfy the HLP axiom.
\end{lem}

Indeed if we have some piece-wise linear curve in $A'=\Pi(A,B)$, we have a
family of piece-wise-linear curves in $B$, starting from $A$. The end-points
give us the necessary homotopy from $B$ to $C$. The theorem is proven. \hfill$\square$

\section{Six-term exact sequences}
Let us recall the definition of the periodic cyclic theory functor
$\HP$, see \cite{connes}, \cite{loday} for more details.
\begin{defn}{\rm
Consider a C*-algebra with unity $1$. We use the maximum norm for
the tensor product of *-algebras. The cyclic complex $C_*(A)$ is
defined as
\begin{equation}C_n(A) := A \otimes \underbrace{\bar{A}
\otimes \dots \otimes \bar{A}}_{\mbox{ n times }}. \end{equation}
Define the Hochschild  operators
\begin{equation}\CD C_0(A) := A @<b,b'<< C_1(A) = A \otimes \bar{A} @<b,b'<< \dots @<b,b'<<
C _n(A) @<b,b'<< \endCD \dots, \end{equation} where $\bar{A} := A
/\mathbf C$. \begin{equation}\b'(a_0 \otimes
a_1\otimes\dots\otimes \otimes a_n) = \sum_{i=0}^{n-1} (-1)^{i}
a_0\otimes \dots\otimes a_ia_{i+1}\otimes\dots\otimes
a_n.\end{equation}
\begin{equation}\b(a_0 \otimes a_1\otimes\dots\otimes
\otimes a_n) = \sum_{i=0}^{n-1} (-1)^{i} a_0\otimes \dots\otimes
a_ia_{i+1}\otimes\dots\otimes a_n + (-1)^n a_na_0\otimes
a_1\otimes\dots\otimes a_{n-1}.\end{equation} It is easy to check
that \begin{equation} \b^2 = 0, \qquad (\b')^2 = 0.\end{equation}
 One define the
cyclic operator $t : C_n(A) \to C_n(A)$ following the formula
\begin{equation}\t(a_0\otimes a_1\otimes\dots \otimes a_n) :=
(-1)^{n+1}a_n\otimes a_0\otimes\dots \otimes a_{n-1}\end{equation}
It is easy to check that \begin{equation} \t^{n+1} = \1
.\end{equation} Define the operators \begin{equation} N = \1 + \t
+ \dots + \t^n.\end{equation} From this cyclic complex we have the
bi-complex $\mathcal C(A)=\{\mathcal C_{p,q} \}$
\begin{equation}\CD
\dots @. \vdots @. \vdots @. \vdots @. \vdots @. \dots \\
@. @VVb'V  @VVbV   @VVb'V  @VVbV @.\\
\dots @>\t-\1>> C_3(A) @>\N>> C_3(A)@>\t-\1>> C_3(A)@>\N>> C_3(A) @>\t-\1>> \dots\\
@. @VVb'V  @VVbV   @VVb'V  @VVbV @.\\
\dots @>\t-\1>> C_2(A) @>\N>> C_2(A)@>\t-\1>> C_2(A)@>\N>> C_2(A) @>\t-\1>> \dots\\
@. @VVb'V  @VVbV   @VVb'V  @VVbV @.\\
\dots @>\t-\1>> C_1(A) @>\N>> C_1(A)@>\t-\1>> C_1(A)@>\N>> C_1(A) @>\t-\1>> \dots\\
@. @VVb'V  @VVbV   @VVb'V  @VVbV @.\\
 \dots @>\t-\1>> C_0(A) @>\N>> C_0(A)@>\t-\1>> C_0(A)@>\N>> C_0(A) @>\t-\1>> \dots
\endCD \end{equation}
The {\it total complex} is defined as
\begin{equation} \Tot\mathcal C(A)_{i} = \prod_{p+q = i (\mod 2)}
 \mathcal C_{p,q},\quad i=0,1\end{equation} with differentials \begin{equation}B = \d_{v}
 + \d_{h}, \end{equation} where $\d_v$ and $\d_h$ are the
 differential following the vertical or horizontal direction,
 correspondingly. The homology of the total complex is called the
 {\it periodic cyclic homology} \begin{equation} \HP_i(A) :=
 \H_i(\Tot\mathcal C(A)), \quad i=0,1. \end{equation}
  }\end{defn}
It is well-known that periodic cyclic homology $\HP_i$ are a
generalized homology functors, introduced by A. Connes
\cite{connes}, \cite{loday}.
\begin{thm}\label{thm2.1}
For every NC Serre fibration $f: A \longrightarrow B$, the corresponding hexagon in periodic cyclic homology holds
\begin{equation}\CD \HP_*(A) @>>> \HP_*(B) @>>> \HP_*(B,A)\\
@A\partial_{*+1} AA @. @VV\partial_* V\\
\HP_{*+1}(B,A) @<<< \HP_{*+1}(B) @<<< \HP_{*+1}(A)
\endCD\end{equation}
\end{thm}
Before prove the theorem we do introduce a noncommutative suspension.
\begin{defn}
Noncommutative (shortly, NC) Suspension $SA$ of $A$ {\rm is by
definition $A \otimes \mathbf C_0(0,1)$.}
\end{defn}
Following the construction of the NC mapping cone of a morphism
$f: A \to B$, we have \cite{diep1}
\begin{lem} There is an exact sequence of algebras
$$\CD \dots @>>> \S^2A @>>> \S\Cone(f:A \to B) @>>> \S\Cyl(f:A\to B)
@>>>  @.\endCD$$
\begin{equation} \CD @. @>>>\S A @>>> \Cone(f:A \to B) @>>> \Cyl(f:A\to B) @>>> A
.\endCD\end{equation}
\end{lem}
The next lemma  is a natural consequence from the well-known
K\"unneth formula an Bott Periodicity of $\HP$.
\begin{lem}
For the periodic cyclic homology $\HP$, there is a natural
isomorphism \begin{equation}\HP_{*}(SA) \cong \HP_{*\pm
1}(A).\end{equation}
\end{lem}

\noindent{\sc Proof of Theorem \ref{thm2.1}}.

\hfill{$\square$}

Because in the category of NCCW-complexes, every morphism is
homotopic to a NC Serre fibration, we have
\begin{cor}
In the category of NCCW-complexes each map $f: A \to B$ admits a
six-term exact sequence
\begin{equation}\CD \HP_*(A) @<<< \HP_*(B) @<<< \HP_*(B,A)\\
@V\partial_{*+1} VV @. @AA\partial_* A\\
\HP_{*+1}(B,A) @>>> \HP_{*+1}(B) @>>> \HP_{*+1}(A).
\endCD\end{equation}
\end{cor}
{\it Proof.} It is well-known Bott Periodicity of type
$\HP_*(\S^2A) \cong \HP_*(A)$ and $\HP_*(\S^2A) \cong \HP_*(A)$.
From the long exact sequence (up-to homotopy) $$\CD  A @<<<
\Cyl(f:A\to B) @<<<  \Cone(i: A\to B)@<<< \S A \endCD$$
\begin{equation}\CD . @<<< \S\Cyl(f:A\to B) @<<< \S\Cone(i: A\to B) @<<<
\S^2A @>>> \dots \endCD\end{equation} we have the first connecting
homomorphism
\begin{equation}\CD \partial_* :  \HP_{*+1}(A)@>>>\HP_*(B,A) \endCD\end{equation}
which gives us the first exact sequence
\begin{equation}\CD \HP_*(A) @<<< \HP_*(B) @<<< \HP_*(B,A)\\
@. @.  @AA\partial_* A \\
@. @. \HP_{*+1}(A). \endCD\end{equation} Then apply for the
suspension $SA$ and use the Bott 2-periodicity in HP-theory we
have the second part
\begin{equation}\CD HP_*(A) \cong HP_{*+2}(A) @. @. \\
@V\partial_{*+1} VV @. @. \\
\HP_{*+1}(B,A) @>>> \HP_{*+1}(B) @>>> \HP_{*+1}(A)
\endCD\end{equation} of the commutative diagram. \hfill$\Box$

For the K-theory we have the analogous results
\begin{thm}\label{thm2.6}
For every NC Serre fibration $f: A \longrightarrow B$, the corresponding hexagon in K-theory holds
\begin{equation}\CD \K_*(A) @<<< \K_*(B) @<<< \K_*(B,A)\\
@V\partial_{*+1} VV @. @AA\partial_* A\\
\K_{*+1}(B,A) @>>> \K_{*+1}(B) @>>> \K_{*+1}(A),
\endCD\end{equation}
\end{thm} and we have also

\begin{cor}
In the category of NC CW-complexes each map $f: A \to B$ admits a six-term exact sequence
\begin{equation}\CD \K_*(A) @<<< \K_*(B) @<<< \K_*(B,A)\\
@V\partial_{*+1} VV @. @AA\partial_* A\\
\K_{*+1}(B,A) @>>> \K_{*+1}(B) @>>> \K_{*+1}(A).
\endCD\end{equation}
\end{cor}

\section{Application}
In the ordinary algebraic topology one used namely the computation of cohomology of spheres to compute homology of orthogonal and unitary groups. In our noncommutative theory, the cyclic homology of quantum orthogonal and quantum unitary groups are known from the representation theory we then use to compute homology of quantum spheres as noncommutative quotients.

\subsection{Quotients of quantum orthogonal groups}
Let us consider the quantum orthogonal groups $C^*_q(\SO(n))$. The natural inclusion $\SO(n-1) \hookrightarrow \SO(n)$ give us the homomorphism
 \begin{equation}C^*_q(\SO(n-1)\setminus\SO(n)) \hookrightarrow  C^*_q(\SO(n)).\end{equation}
 We have therefore up to homotopy a NC Serre fibration what for simplicity we denote
 again $(C^*_q(\SU(n):C^*_q(\SU(n-1))):= $
 \begin{equation}\Pi(C^*_q(\SO(n-1)\setminus\SO(n),\Cyl(,C^*_q(\SO(n)),C^*_q(\SO(n-1)\setminus\SO(n)))\end{equation}
 and our associated exact sequence (up to homotopy of algebras) as
 {\scriptsize
 $$ \CD C^*_q(\SO(n-1)\setminus\SO(n)) @<<<  C^*_q(\SO(n)) @<<< \Cone(C^*_q(\SO(n-1)\setminus\SO(n)),C^*_q(\SO(n)))@<<<
 @.\endCD$$
 $$\CD \S C^*_q(\SO(n-1)\setminus\SO(n)) @<<<  \S C^*_q(\SO(n)) @<<<
 \S\Cone(C^*_q(\SO(n-1)\setminus\SO(n)),C^*_q(\SO(n)))
 \endCD$$
 \begin{equation}\CD @. @<<< \S^2C^*_q(\SO(n-1)\setminus\SO(n))@<<< \dots \endCD
 \end{equation}}
We can now apply the six-term exact sequence for Serre fibrations
and have {\scriptsize
\begin{equation}\CD \HP_*(C^*_q(\SO(n-1)\setminus\SO(n))) @<<< \HP_*(C^*_q(\SO(n)))
@<<< \HP_*(C^*_q(\SO(n)):C^*_q(\SO(n-1)))\\
@V\partial_{*+1} VV @. @AA\partial_* A\\
\HP_{*+1}(C^*_q(\SO(n)):C^*_q(\SO(n-1))) @>>>
\HP_{*+1}(C^*_q(\SO(n))) @>>>
\HP_{*+1}(C^*_q(\SO(n-1)\setminus\SO(n)))\endCD\end{equation}} In
this six-term exact sequences the $\HP_*(C^*_q(\SO(n))) \cong
H^*_{DR}(\mathbb T)^W$ are well-known, where $\mathbb T$ is a
fixed maximal torus in $\SO(n)$ and $W = W(\mathbb T) = \mathcal
N(\mathbb T)/\mathbb T$ is the Weyl group of this maximal torus,
see \cite{diepkukutho1}-\cite{diepkukutho2}. We can fix an
immersion $\SO(n-1) \hookrightarrow \SO(n)$ so that the
corresponding tori are included by immersions one-into-another.

\subsection{Quotients of quantum unitary groups}
By analogy, we  consider the quantum unitary groups $C^*_q(\SU(n))$. The natural inclusion $\SU(n-1) \hookrightarrow \SU(n)$ give us the homomorphism
 \begin{equation}C^*_q(\SU(n-1)\setminus\SU(n)) \hookrightarrow  C^*_q(\SU(n)).\end{equation}
 We have therefore up to homotopy a NC Serre fibration what for simplicity we denote
 again $(C^*_q(\SU(n):C^*_q(\SU(n-1))):=$
 \begin{equation}\Pi(C^*_q(\SU(n-1)\setminus\SU(n),\Cyl(,C^*_q(\SU(n),C^*_q(\SU(n-1)\setminus\SU(n)))\end{equation}
 and our associated exact sequence (up to homotopy of algebras) as
 {\scriptsize
$$ \CD C^*_q(\SU(n-1)\setminus\SU(n)) @<<<  C^*_q(\SU(n)) @<<< \Cone(C^*_q(\SU(n-1)\setminus\SU(n)),C^*_q(\SU(n)))@<<<
 @.\endCD$$
 $$\CD \S C^*_q(\SU(n-1)\setminus\SU(n)) @<<<  \S C^*_q(\SU(n)) @<<<
 \S\Cone(C^*_q(\SU(n-1)\setminus\SU(n)),C^*_q(\SU(n)))
 \endCD$$
 \begin{equation}\CD @. @<<< \S^2C^*_q(\SU(n-1)\setminus\SU(n))@<<< \dots \endCD
 \end{equation}}

We can now also apply the six-term exact sequence for Serre
fibrations and have {\scriptsize
\begin{equation}\CD \HP_*(C^*_q(\SU(n-1)\setminus\SU(n))) @<<< \HP_*(C^*_q(\SU(n))) @<<< \HP_*(C^*_q(\SU(n)):C^*_q(\SU(n-1)))\\\
@V\partial_{*+1} VV @. @AA\partial_* A\\
\HP_{*+1}(C^*_q(\SU(n)):C^*_q(\SU(n-1))) @>>>
\HP_{*+1}(C^*_q(\SU(n))) @>>>
\HP_{*+1}(C^*_q(\SU(n-1)\setminus\SU(n)))\endCD\end{equation}} In
this six-term exact sequences the $\HP_*(C^*_q(\SU(n))) \cong
H^*_{DR}(\mathbb T)^W$ are well-known, where $\mathbb T$ is a
fixed maximal torus in $\SU(n)$ and $W = W(\mathbb T) = \mathcal
N(\mathbb T)/\mathbb T$ is the Weyl group of this maximal torus,
see \cite{diepkukutho1}-\cite{diepkukutho2}. We can fix an
immersion $\SU(n-1) \hookrightarrow \SU(n)$ so that the
corresponding tori are included by immersions one-into-another.

\section*{Acknowledgments}
The work was supported in part by Vietnam National Project for
Research in Fundamental Sciences and was completed  during the
stay in June and July, 2007 of the author, in Abdus Salam ICTP,
Trieste, Italy.  The author expresses his deep and sincere thanks
to Abdus Salam ICTP and especially Professor Dr. Le Dung Trang for
the invitation and for providing the nice conditions of work, and
Professor C. Schochet for some Email discussion messages and in
particular for the reference \cite{schochet}.

{\sc Institute of Mathematics, National Centre for Science and
Technology of Vietnam, 18 Hoang Quoc Viet Road, Cau Giay district
, 10307 Hanoi, Vietnam}
\\
{\tt Email: dndiep@math.ac.vn}
\end{document}